\documentclass[10pt]{article}
\usepackage{epsfig}

\def\R{\hbox{\bf R}}

\def\N{\hbox{\bf N}}

\def\dist{\hbox{\rm dist}}

\def\varphi0{{\cal Q}}

\def \proof{{\sl Proof: }}

\newcommand{\be}{\begin{equation}}
\newcommand{\ee}{\end{equation}}
\newcommand{\baa}{\begin{array}}
\newcommand{\eaa}{\end{array}}
\newcommand{\ba}{\begin{eqnarray}}
\newcommand{\ea}{\end{eqnarray}}

\newtheorem{theo}{\textsc{Theorem}}[section]
\newtheorem{lem}[theo]{\textsc{Lemma}}
\newtheorem{pro}[theo]{\textsc{Proposition}}
\newtheorem{cor}[theo]{\textsc{Corollary}}
\newtheorem{defi}[theo]{\textsc{Definition}}
\newtheorem{rem}[theo]{\textsc{Remark}}

\begin{document}
\title{{\bf \Large An Unstable Elliptic Free Boundary Problem\\ arising in Solid Combustion}}
\author{
\normalsize\textsc{R. Monneau}\\
{\normalsize\it Ecole Nationale des Ponts et Chauss\'ees, CERMICS,
 6 et 8 avenue
  Blaise Pascal,}\\
{\normalsize\it  Cit\'e Descartes Champs-sur-Marne, 77455 Marne-la-Vall\'ee  Cedex 2, France}\\
\normalsize\textsc{G.S. Weiss}\\
{\normalsize\it Graduate School of Mathematical Sciences,}\\
{\normalsize\it
University of Tokyo,
3-8-1 Komaba, Meguro, Tokyo,
153-8914 Japan,}\\
{\normalsize\it Guest of the Max Planck Institute
 for Mathematics in the Sciences,}\\
{\normalsize\it Inselstr. 22, D-04103 Leipzig, Germany}
\thanks{G.S. Weiss has been partially supported by the Grant-in-Aid
15740100 of the Japanese Ministry of Education and partially supported
by a fellowship of the Max Planck Society. Both authors thank the Max Planck
Institute for Mathematics in the Sciences for the hospitality
during their stay in Leipzig.}\\
}
%\date{}
\maketitle

%%%%%%%%%%%%%%%%%%%%%%%%%%%%%%%%%%%%%%%%%%%%%%%%%%%%%%%%%%%%%%%%%%%%%%%%%%%%%%%%
%%%%%%%%%%%%%%%%%%%%%%%%%%%%%%%%%%%%%%%%%%%%%%%%%%%%%%%%%%%%%%%%%%%%%%%%%%%%%%%%

\centerline{\small{\bf{Abstract}}}
We prove a regularity result for the
unstable elliptic
free boundary problem
\begin{equation}
\Delta u = -\chi_{\{ u>0\}}
\end{equation}
related to traveling waves in a problem arising in solid combustion.
\\
The maximal solution and every local minimizer of the
energy are regular, that is, $\{ u=0\}$ is locally 
an analytic surface and
$u|_{\overline{\{ u>0\}}}, u|_{\overline{\{ u<0\}}}$ are locally analytic functions.
Moreover we prove a partial regularity result for 
solutions that are non-degenerate of second order:
here $\{ u=0\}$ is analytic up to a closed set of Hausdorff
dimension $n-2$. We discuss possible singularities.\\
\noindent{\small{}}\hfill\break

\noindent{\small{\bf{AMS Classification:}}} {\small{35R35, 35J60, 35B65.}}\hfill\break
\noindent{\small{\bf{Keywords:}}} {\small{
free boundary, regularity, monotonicity formula, frequency, solid combustion, 
singularity, unstable problem, Aleksandrov reflection, unique blow-up limit,
second variation, maximal solution}}\hfill\break

%%%%%%%%%%%%%%%%%%%%%%%%%%%%%%%%%%%%%%%%%%%%%%%%%%%%%%%%%%%%%%%%%%%%%%%%%%%%%%%%
%%%%%%%%%%%%%%%%%%%%%%%%%%%%%%%%%%%%%%%%%%%%%%%%%%%%%%%%%%%%%%%%%%%%%%%%%%%%%%%%

\section{Introduction}
We prove a regularity result for the unstable elliptic
free boundary problem
\begin{equation}\label{eq}
\Delta u = -\chi_{\{ u>0\}}\; .
\end{equation}
The problem (\ref{eq}) is related to traveling wave solutions in solid
combustion with ignition temperature:
let us consider the solid combustion system
\begin{equation}
\begin{array}{l}
\partial_t \theta - \Delta \theta = g(\eta)\chi_{\{\theta>\theta_0\}},\\
\partial_t \eta = g(\eta)\chi_{\{\theta>\theta_0\}}
\end{array}
\end{equation}
where $\theta$ is the temperature,
$\theta_0$ is the ignition temperature and $\eta \in (0,1)$ is
the fractional conversion.
Although an ignition temperature has no meaning for
gas flames, it has been recently rediscovered and
used in combustion synthesis (see for example \cite{varma},\cite{beck},\cite{khaikin} and \cite{solomonov}).
The reaction kinetics suggested in \cite[p.1462]{varma}
is $g(\eta)=k_0\> \chi_{\{\eta<1\}}$, but the particular form
of $g$ is not of importance for what follows.
Solving the ODE
we obtain for $u := \theta-\theta_0$
\begin{equation}\label{sing}
\partial_t u -\Delta u  = c\chi_{\{ u>0\}}\; ,
\end{equation}
where $c$ is a non-negative memory term depending on $(t,x)$ as well
as the history $\{ u(s,x), s<t\}.$
Although $c$ is important when considering the large-time
behavior of $u,$ we may consider (\ref{sing}) to be a 
perturbation of
\begin{equation}\label{parab}
\partial_t u -\Delta u  = \chi_{\{ u>0\}}
\end{equation}
when we are interested in transient or local phenomena
at the ignition front $\partial \{u<0\}$.
Actually the traveling pulses in our model correspond
well to the fingering phenomenon for burned regions
observed in solid combustion experiments (see for example \cite{zik}). 
\\
We are interested in traveling wave solutions. As we deal
in the present paper mostly with regularity issues, we
may drop the drift term resulting from the time derivative
of the traveling wave.
\\
An equation similar to our elliptic one arises
in the composite membrane problem (see \cite{chanillo1},
\cite{chanillo2}, \cite{blank}).
Another application is the shape of self-gravitating rotating
fluids describing stars (see \cite[equation (1.26)]{stars}). 
\\
From a mathematical point of view, (\ref{parab})
is the equation of the parabolic obstacle problem
with inverted sign, and (\ref{eq})
is the equation of the elliptic obstacle problem
with inverted sign. The change of sign changes the
character of the problem {\em drastically} in that it changes
the {\em stable} obstacle problem into an {\em unstable
problem.} In (\ref{eq}) and (\ref{parab}), we find
examples of non-uniqueness, bifurcation phenomena etc.
\\
As surprisingly many known free boundary problems turn
out to be stable problems, this means, unfortunately,
that many known methods in free boundary problems
do not apply here. Examples of PDE techniques that do not
work are, apart from all one-phase methods, the
Bernstein technique, the
Alt-Caffarelli-Friedman monotonicity formula (\cite{ACF}) and
the differential inequality technique of Cazenave-Lions
(\cite{cazenavelions}).
\\
Our main result is 
that the maximal solution and every local minimizer of the
energy are regular, that is, $\{ u=0\}$ is locally an
analytic surface and
$u|_{\overline{\{ u>0\}}}, u|_{\overline{\{ u<0\}}}$ are
locally analytic functions (Theorem \ref{main}). 
\\
The surprise is that -- in contrast to the usual
procedure -- we obtain $C^{1,1}$-regularity of local
minimizers {\em by proving regularity of the free boundary
first!}\\
For general solutions that are non-degenerate of second order, we prove
a partial regularity result:
here $\{ u=0\}$ is smooth up to a set of Hausdorff
dimension $n-2$ (Proposition \ref{partial}).
We discuss the behavior at possible singularities.\\
In case of a non-degenerate minimal solution
in two dimensions we also obtain that $\{ u=0\}$ 
consists of Lipschitz arcs meeting at right angles
in at most finitely many singularities.
\\
{\bf Acknowledgment:} We thank Ivan Blank, Carlos Kenig, Herbert Koch
and Stephan Luckhaus
for fruitful discussions. 
\section{Notation}
Throughout this article $\R^n$ will be equipped with the Euclidean
inner product $x\cdot y$ and the induced norm $\vert x \vert\> .$
We define $e_i$ as the $i$-th unit vector in $\R^n\> ,$ and
$B_r(x^0)$ will denote the open $n$-dimensional ball of center
$x^0\> ,$ radius $r$ and volume $r^n\> \omega_n\> .$
We shall often use
abbreviations for inverse images like $\{u>0\} := 
\{x\in \Omega\> : \> u(x)>0\}\> , \> \{x_n>0\} := 
\{x \in \R^n \> : \> x_n > 0\}$ etc. 
and occasionally 
we shall employ the decomposition $x=(x_1,\dots,x_n)$ of a vector $x\in \R^n\> .$
We will use the $k$-dimensional Hausdorff measure
${\cal H}^k$ approximated by ${\cal H}^{k,\delta}$ which we define
as the $H^\delta_k$ of \cite{giu}.
When considering a set $A\> ,$ $\chi_A$ shall stand for
the characteristic function of $A\> ,$ 
while
$\nu$ shall typically denote the outward
normal to a given boundary.
\section{Existence and Non-Degeneracy} 
\begin{lem}{\bf (Existence of a maximal and a minimal solution)}\\
Let $\Omega$ be a bounded domain of class $C^{2,\alpha}$ in
$\R^n,$ and assume that the Dirichlet
boundary data $u_D\in C^{2,\alpha}(\bar\Omega).$\\
Then there exists a {\em maximal solution} $u$ with
the following properties:
$u\in W^{2,p}(\Omega)$ for every $p\in [1,+\infty),$
$u=u_D$ on $\partial \Omega$ and $u\ge v$
for every subsolution $v\in W^{2,n}(\Omega)$ of (\ref{eq}) in $\Omega'\subset\Omega$
satisfying $v\le u$ on $\partial \Omega'.$\\
There also exists a {\em minimal solution} with analogous
properties.
\end{lem}
\proof
We prove the existence of the maximal solution:
We consider a regularization of the equation from above,
$$ \Delta v = -\bar \beta_\epsilon(v),$$
where $\bar \beta_\epsilon\in C^\infty(\R),\bar \beta_\epsilon(z)\ge 
\chi_{\{ z>0\}}=: \beta(z)$ in $\R$ and $\bar\beta_\epsilon
\downarrow \beta$ as $\epsilon\downarrow 0.$
By Perron's method, there exists a maximal solution $u_\epsilon$ in the above
sense for each $\epsilon>0.$
By $W^{2,p}$-estimates the family $(u_\epsilon)_{\epsilon\in (0,1)}$
is bounded in $W^{2,p}(\Omega).$
Moreover, for any subsolution $v$ of (\ref{eq}) and
$0<\tilde \epsilon<\epsilon$ we obtain
$$\Delta v \geq - \beta(v) \geq -\bar\beta_{\epsilon}(v)\hbox{ and }
\Delta u_{\tilde \epsilon} = - \beta_{\tilde \epsilon}(u_{\tilde \epsilon}) 
\geq -\bar\beta_{\epsilon}(u_{\tilde \epsilon})$$
so that $v$ and $u_{\tilde \epsilon}$ are subsolutions of the $\epsilon$-equation.
Consequently $v\le u_\epsilon$ and $u_{\tilde \epsilon}\le u_\epsilon.$
Using the fact that $D^2u=0$ a.e. on $\{ u=0\}$,
it follows that $u(x) := \lim_{\epsilon\to 0} u_\epsilon(x)$ is the
maximal solution.
\begin{lem}\label{maxndeg}
There exists a positive constant $c_n$ depending only on the
dimension $n,$ such that the maximal solution $u$ with respect
to given boundary data satisfies
$$ B_r(x^0)\subset \Omega, \inf_{\partial B_r(x^0)} u > -c_n r^2 \Rightarrow
u(x^0)>0 \; .$$
Thus for any $x^0\in \{ u=0\}$ and $r<\dist(x^0,\Omega^c)$
we obtain $\inf_{\partial B_r(x^0)} u \le -c_n r^2.$
\end{lem}
\proof
Compare $u(x^0+rx)/r^2$ for suitable $\theta$ to the ``stationary pulse''
$p(\theta x)/\theta^2,$ where
$$p(x) = \left\{\begin{array}{ll}
(1-\vert x \vert^2)/4, & \vert x \vert \le 1,\\
-\log(\vert x \vert^2)/4, & \vert x \vert > 1\end{array}\right.$$
in the case $n=2$ and to
$$p(x) = \left\{\begin{array}{ll}
{1\over {2n}}(1-\vert x \vert^2), & \vert x \vert \le 1,\\
{1\over {n(2-n)}}(\vert x \vert^{2-n}-1), & \vert x \vert > 1\end{array}\right.$$
in the case $n>2.$\\
\begin{rem}For the case of two dimensions the optimal constant $c_2=1/(4e).$
\end{rem}
\begin{defi}[Non-degeneracy]
Let $u$ be a solution of (\ref{eq}) in $\Omega,$
satisfying at $x^0\in \Omega$
\begin{equation}\label{ndeg}
\liminf_{r\to 0} r^{-2}\left(\int_{\partial B_{r_m}(x^0)} u^2\> d{\cal H}^{n-1}
\right)^{1\over 2}>0\; .
\end{equation}
Then we call $u$ ``non-degenerate of second order at $x^0$''.
We call $u$ ``non-degenerate of second order'' if it is
non-degenerate of second order at each point in $\Omega$. 
\end{defi}
\begin{rem}
The maximal solution is non-degenerate of second order.
\end{rem}
\begin{lem}\label{minndeg}
Each minimizer $u$ of the energy
\[ E(v) = \int_{B_{r_0}(x^0)} |\nabla v|^2 - 2\max(v,0)\]
in $K:= \{ v\in W^{1,2}(B_{r_0}(x^0)) : v=u \hbox{ on } \partial
B_{r_0}(x^0)\}$
satisfies at $x^0$ the second-order non-degeneracy property (\ref{ndeg}).
\end{lem}
{\em Proof:}
Define for $r \in (0,r_0)$ the solution $v := u(x^0+rx)/r^2$ and let $p$
be the stationary pulse in $B_1(0)$
with boundary data $\inf_{\partial B_1(0)} v$. Comparing the energy
of $v$ to that of $w:= \max(p,v)$ we obtain
\\
\[ 0\ge \int_{B_1(0)}  |\nabla v|^2 - |\nabla w|^2-2\max(v,0) 
+ 2\max(w,0)\]\[
=
\int_{B_1(0)} -(\Delta v + \Delta w)(v-w) + 2(\max(w,0)-\max(v,0))
\]\[
=
\int_{B_1(0)\cap\{ v\le0\} \cap \{ p>0\}} v+p\; .\]
Now assume that
$v_m$ is a sequence of minimizers such that
$v_m(0)=0$ and
$\inf_{\partial B_1(0)} v_m \to 0$ as $m\to\infty$.
Let $\epsilon\in (0,1)$ be fixed.
In the case that ${\cal L}^n(B_{1-\epsilon}(0)\cap
\{ v_m\le 0\})\not\to 0$ as $m\to\infty$,
we obtain immediately a contradiction since
$p_m$ is for large $m$ on 
$B_{1-\epsilon}(0)$ uniformly estimated from below by a positive constant
depending only on $n$ and $\epsilon$.
In the case ${\cal L}^n(B_{1-\epsilon}(0)\cap
\{ v_m\le 0\})\to 0$ as $m\to\infty$,
we know by $\lim_{m\to\infty}\inf_{B_1(0)}v_m \to 0$
and $\int_{\partial B_s(0)} v_m \> d{\cal H}^{n-1}\le 0$ for every $s\in (0,1]$ (which follows from $v_m$ being superharmonic and $0$ at the origin)
that $\int_{B_1(0)} |v_m|
\to 0$ as $m\to\infty$.
Consequently
\[ 0 \gets \int_{B_{1-\epsilon}(0)} \Delta v_m = \int_{B_{1-\epsilon}(0)} -\chi_{\{ v_m>0\}}
\le -c_1 <0\; \hbox{as } m\to\infty,\]
a contradiction. Thus $u(x^0)=0$ implies
\[ \inf_{B_1(0)} {u(x^0+r\cdot )\over {r^2}} \le -c_2 <0\]
for all $r\in (0,r_0)$.
\begin{lem}\label{maxmin}
The maximal solution $u_{\rm max}$ is {\em the} maximal minimizer of the
energy
\[ E(v) = \int_{\Omega} |\nabla v|^2 - 2\max(v,0)\]
in $K:= \{ v\in W^{1,2}(\Omega) : v=u_D \hbox{ on } \partial
\Omega\}$.
\end{lem} 
{\em Proof:}
For any $v\in K$ we get
\[ 0\ge \int_\Omega \max(v,0)-\max(u_{\rm max},0)+u_{\rm max}\chi_{\{v>0\}}
- v\chi_{\{ u_{\rm max}>0\}}\]\[
= \int_\Omega (\chi_{\{v>0\}}+\chi_{\{ u_{\rm max}>0\}})(u_{\rm max}-v)
+ 2(\max(v,0)- \max(u_{\rm max},0))\]\[
= \int_\Omega |\nabla u_{\rm max}|^2 - |\nabla v|^2 + 2\max(v,0)- 2\max(u_{\rm max},0)
\; .\]
\section{Monotonicity Formula and Frequency Lemma}
A powerful tool is now a monotonicity formula introduced 
in \cite{cpde} by one of the authors for a class of semilinear free boundary
problems. For the sake of completeness let us state 
the unstable case here:
\begin{theo}[Monotonicity formula]
\label{mon}
Suppose that $u$ is a solution of (\ref{eq}) in $\Omega$
and that $B_\delta(x^0)\subset \Omega\> .$
Then for all $0<\rho<\sigma<\delta$
the function 
\[ \Phi_{x^0}(r) := r^{-n-2} \int_{B_r(x^0)} \left( 
{\vert \nabla u \vert}^2 \> -\> 2\max(u,0)
\right)\]\[ 
- \; 2 \> r^{-n-3}\>  \int_{\partial B_r(x^0)}
u^2 \> d{\cal H}^{n-1}\; ,\]
defined in $(0,\delta)\> ,$ satisfies the monotonicity formula
\[ \Phi_{x^0}(\sigma)\> -\> \Phi_{x^0}(\rho) \; = \;
\int_\rho^\sigma r^{-n-2}\;
\int_{\partial B_r(x^0)} 2 \left(\nabla u \cdot \nu - 2 \>
{u \over r}\right)^2 \; d{\cal H}^{n-1} \> dr \; \ge 0 \; \; .\] 
\end{theo}
We also need the following frequency lemma, which has been
proven in \cite[Lemma 4.1]{interfaces}:
\begin{lem}[Frequency lemma]\label{freq}
Let $\alpha-1 \in \N\> ,$ let $
w\in W^{1,2}(B_1(0))$ be a harmonic function in $B_1(0)$
and assume that $D^j w(0)=0$ for $0\le j \le \alpha - 1\> .$\\
\[ \hbox{Then } \int_{B_1(0)} {\vert \nabla w \vert}^2 
\> - \> \alpha \int_{\partial B_1(0)} w^2 \> d{\cal H}^{n-1}
\; \ge \; 0\; ,\]
and equality implies that $w$ is homogeneous of degree
$\alpha$ in $B_1(0)\> .$
\end{lem}
\section{Classification of Blow-up Limits}
A result related to the following 
classification of blow-up limits
is contained in \cite[Theorem 2.5]{blank}.
Note however that the proofs of the related parts
are largely different.
\begin{pro}[Classification of blow-up limits with fixed center]\label{fixedcenter}
Let $u$ be a solution of (\ref{eq}) in $\Omega$
and let us consider a point
$x^0\in \Omega\cap \{ u=0\}\cap\{ \nabla u =0\}.$\\
1) In the case $\Phi_{x^0}(0+)=-\infty$,
$\lim_{r\to 0} r^{-3-n}\int_{\partial B_r(x^0)} u^2 \> d{\cal H}^{n-1}
= +\infty$, and for
$S(x^0,r) := \left(r^{1-n}\int_{\partial B_{r}(x^0)} u^2\> d{\cal H}^{n-1}
\right)^{1\over 2}$
each limit of
\[ \frac{u(x^0+r x)}{S(x^0,r)}\]
as $r\to 0$ is a homogeneous harmonic polynomial of degree $2$.
\\
2) In the case $\Phi_{x^0}(0+)\in (-\infty,0)$,
\[ u_r(x) := \frac{u(x^0+r x)}{r^2}\]
is bounded in $W^{1,2}(B_1(0))$,
and each limit as $r\to 0$ is a homogeneous solution of degree $2$.\\
3) Else $\Phi_{x^0}(0+)=0$, and
\[ \frac{u(x^0+r x)}{r^2}\to 0\hbox{ in } W^{1,2}(B_1(0)) \hbox{ as } r\to 0\; .\]
\end{pro}
\proof
In all three cases,
\[ \partial_r {1\over 2}\int_{\partial B_1(0)} u_r^2 \> d{\cal H}^{n-1}
\; = \; \int_{\partial B_1(0)} u_r \partial_r u_r\> d{\cal H}^{n-1} 
\]\[ = \; {1\over r}  \int_{\partial B_1(0)} u_r (\nabla u_r\cdot x - 2u_r)\> d{\cal H}^{n-1} 
\]\[ = \; {1\over r} \left(\int_{B_1(0)} -\max(u_r,0) + |\nabla u_r|^2 \; - \; 2
\int_{\partial B_1(0)} u_r^2 \> d{\cal H}^{n-1}\right)
\]\[ = \; {1\over r} \left(\Phi_{x^0}(r) \; + \; \int_{B_1(0)} \max(u_r,0)\right)\; .\]
In particular, 
\begin{equation}\label{blowup1}
\int_{B_1(0)} \max(u_r,0)\ge M := - \Phi_{x^0}(0+)
\hbox{ implies }\partial_r \int_{\partial B_1(0)} u_r^2 \> d{\cal H}^{n-1}\ge 0\; .
\end{equation}
Observe now that 
\begin{equation}\label{blowup2}
\int_{\partial B_1(0)} \max(-u_{r},0)^2 \> d{\cal H}^{n-1}
\le C_1 \left(1+\int_{\partial B_1(0)} \max(u_{r},0)^2 \> d{\cal H}^{n-1}\right)\; :
\end{equation}
supposing towards a contradiction that this is not true,
$\int_{\partial B_1(0)} u_r^2 \> d{\cal H}^{n-1}$
must be unbounded for a sequence $(r_m)_{m\in \N}$
not satisfying (\ref{blowup2}).
We may divide the function of the monotonicity formula \ref{mon}
\[ \int_{B_1(0)} |\nabla u_{r_m}|^2 - 2\int_{\partial B_1(0)} u_{r_m}^2 \> d{\cal H}^{n-1}  
\le \Phi_{x^0}(1) \; + \; 2 \int_{B_1(0)}\max(u_{r_m},0)\]
by $\int_{\partial B_1(0)} \max(-u_{r_m},0)^2 \> d{\cal H}^{n-1}$.
As a subsequence $r_m\to 0$, we obtain a weak $L^2(\partial B_1(0))$-limit $v$ of 
$v_m:=u_{r_m}/\left(\int_{\partial B_1(0)} \max(-u_{r_m},0)^2 \> d{\cal H}^{n-1}\right)^{1/2}$
such that by $\int_{B_1(0)} |u_{r_m}| \le C_2 
\left(\int_{\partial B_1(0)} |u_{r_m}|^2\> d{\cal H}^{n-1}\right)^{1\over 2}$
(which follows from $\min(u_{r_m},0)^2$ being subharmonic and $u_{r_m}$ being
superharmonic and $0$ at the origin) 
\begin{equation}\label{blowup5}
\int_{B_1(0)} |\nabla v_m|^2 - 2\int_{\partial B_1(0)} v_m^2 \> d{\cal H}^{n-1} 
\; \le \; o(1) \hbox{ as } m\to\infty\; ,\end{equation}
and $v$ is in $B_1(0)$ a non-positive harmonic function satisfying $v(0)=0$,
implying by the strong maximum principle that $v=0$ in $B_1(0)$.
That poses a contradiction to the (by (\ref{blowup5})) strong $L^2(\partial B_1(0))$-convergence
of $v_m$ and the fact that
\[ \liminf_{m\to\infty}\int_{\partial B_1(0)} v_m^2\> d{\cal H}^{n-1} >0\; .\]
Note also that in the case $\Phi_{x^0}(0+)>-\infty$, for small $r$
and $\tilde r\in (r/2,r)$
\begin{equation}\label{blowup3}
\begin{array}{l}
1>\Phi_{x^0}(r)-\Phi_{x^0}(\tilde r)
= \int_{\tilde r}^r 2s \int_{\partial B_1(0)} (\partial_s u_s)^2 \> d{\cal H}^{n-1}\> ds
\\
\ge \int_{\partial B_1(0)} (u_r - u_{\tilde r})^2 d{\cal H}^{n-1}\; .\end{array}
\end{equation}
Combining (\ref{blowup3}), (\ref{blowup2}) and (\ref{blowup1})
we see that in the case $\Phi_{x^0}(0+)>-\infty$,
\begin{equation}\label{blowup4}
\int_{\partial B_1(0)} u_r^2 \ge \tilde M
\hbox{ implies }\partial_r \int_{\partial B_1(0)} u_r^2 \> d{\cal H}^{n-1}\ge 0\; .
\end{equation}
But then $\int_{\partial B_1(0)} u_r^2 \> d{\cal H}^{n-1}$ has to be bounded
in the case $\Phi_{x^0}(0+)>-\infty$.
\\
From the monotonicity formula Theorem \ref{mon} we infer that
in the case $\Phi_{x^0}(0+)>-\infty$,
each limit $u_0$ of $u_r$ in $B_1(0)$ is a homogeneous solution 
of degree $2$.\\
Moreover, in the case $\Phi_{x^0}(0+)\ge 0$, 
we obtain that each limit $u_0$ satisfies
\[ 0\le \Phi_{x^0}(0+)
\; = \; \int_{B_1(0)} |\nabla u_0|^2 -2\max(u_0,0)
\> - \> 2 \int_{\partial B_1(0)} u_0^2 \> d{\cal H}^{n-1}
\; = \; -\int_{B_1(0)} \chi_{\{ u_0>0\}}\; .\]
It follows that $u_0\equiv 0$ in $B_1(0)$ and that
$\Phi_{x^0}(0+)=0$.
\\
Last, in the case $\Phi_{x^0}(0+)=-\infty$ we obtain 
that
\[ \lim_{r\to 0} S(x^0,r)/r^2 = \lim_{r\to 0} \left(r^{-3-n}\int_{\partial B_r(x^0)} u^2 \> d{\cal H}^{n-1}\right)^{1\over 2}
= +\infty\; .\]
Taking a subsequence 
\[ w_m(x) := \frac{u(x^0+r_m x)}{S(x^0,r_m)}\]
that converges weakly in $L^2(\partial B_1(0))$ 
to $w_0$ and setting $T_m := S(x^0,r_m)/{r_m}^2$,
we infer from
the monotonicity formula Theorem 
\ref{mon} for large $m$ that
$$\int_{B_1(0)} \vert \nabla u_{r_m}\vert^2
\leq \Phi_{x^0}(r_0) \> + \> \int_{B_1(0)} 2\max(u_{r_m},0)
\> + \> 2 \int_{\partial B_1(0)} u_{r_m}^2 \> d{\cal H}^{n-1}\; .$$
Division by $T_m^2$ yields
$$\int_{B_1(0)} \vert \nabla w_m\vert^2
\leq T_m^{-2}\Phi_{x^0}(r_0) \> + \> T_m^{-1}\int_{B_1(0)} 2\max(w_m,0)
\> + \> 2 \int_{\partial B_1(0)} w_m^2 \> d{\cal H}^{n-1}\; .$$
Since $|\Delta w_m|\le 1$ in $B_{1}(0)$, it follows that 
\begin{equation}
\int_{B_1(0)} \vert \nabla w_0\vert^2
\leq 
2 \int_{\partial B_1(0)} w_0^2 \> d{\cal H}^{n-1}\; ,\end{equation}
that $w_0(0)=\vert \nabla w_0(0)\vert = 0$ 
and that $w_0$ is harmonic in $B_1(0)$.
From Lemma \ref{freq} we infer that
$w_0$ is a homogeneous harmonic polynomial
of degree $2.$
\begin{lem}\label{nohomog}
In two dimensions, the only solution of (\ref{eq})
that is homogeneous of degree $2,$ is the trivial solution
$0.$
\end{lem}
\proof
Let $u$ be a solution of (\ref{eq}) on $\R^2$
that is homogeneous of degree $2.$
Passing to the ODE
$y'' +4y = -\chi_{\{y>0\}},$ each component of $\{ u<0\}$ must
be a cone of opening $\pi/2$ and each component of $\{ u>0\}$ must
for some $\tau\in (0,+\infty)$
be a cone of opening $|\hbox{arcsin}(1/(4\tau))|<\pi/2,$ a contradiction.
\begin{rem}
Let $u$ be a solution in a neighborhood of $x^0\in \Omega$.
If the point $x^0$ is non-degenerate of second order, then
all points in some open neighborhood
of $x^0$ are non-degenerate of second order.
\end{rem}
\proof
$\Phi_{x^0}(0+)<0$ implies by upper semicontinuity 
of \[ x\mapsto\Phi_{x}(0+) \]
that every point $x\in \{ u=0\}$ in some open neighborhood of $x^0$
satisfies 
$\Phi_{x}(0+)<0$ and is therefore non-degenerate of second order. 
\section{Partial Regularity}
A result related to the following Corollary has been independently
obtained in \cite[Theorem 1.1]{henrik}.
\begin{cor}[Partial regularity in two dimensions]\label{isolated}
Let $n=2$ and let $u$ be a solution of (\ref{eq}) in $\Omega$ 
that is non-degenerate of second order. \
Then for each $K\subset\subset \Omega,$
the singular set $K \cap \{ u=0\}\cap
\{\nabla u=0\}$ contains at most finitely many points.
\end{cor}
\proof
Suppose this is not true. Then there is a sequence
$\Omega\cap \{ u=0\} \cap
\{\nabla u=0\}\ni x^m\to x^0\in \Omega\cap \{ u=0\} \cap
\{\nabla u=0\}.$
Take a blow-up limit $u_0$ with respect to the fixed center $x^0$
such that $\partial B_1(0)$ contains a point of
$\{ u_0=0\} \cap
\{\nabla u_0=0\}.$ By Proposition \ref{fixedcenter}
and Lemma \ref{nohomog} 
we know that $u_0$ is a homogeneous harmonic polynomial
of degree $2.$ This is a contradiction, since for a homogeneous harmonic polynomial
of degree $2$ in two dimensions the set $\{ u_0=0\} \cap
\{\nabla u_0=0\}=\{ 0\}.$
\begin{lem}\label{dist}
Let $u$ be a solution of (\ref{eq}) in $\Omega$
that is non-degenerate of second order, let
$x^0\in \Omega\cap \{ u=0\}\cap\{ \nabla u=0\},$
and let $u_0$ be a blow-up limit of 
$$u_m(x) := \frac{u(x^0+r_m x)}{S(x^0,r_m)}$$
in sense of Proposition \ref{fixedcenter}.
Then for each compact set $K\subset \R^n$ and each open set $U\supset
K\cap S_0$ there exists $m_0<\infty$ such that $S_m \cap K \subset U$
for $m\ge m_0\> ;$ here $S_0 := 
\{ u_0=0\}\cap \{\nabla u_0 = 0\}$
and $S_m := 
\{ u_m=0\}\cap \{\nabla u_m = 0\}\> .$ 
\end{lem}
\proof
Suppose towards a contradiction that $S_m \cap (K-U)\ni x^m\to \bar x$
as $m \to \infty \> .$ Then $\bar x\in \{u_0=0\}
\cap \{\nabla u_0=0\}\cap (K-U)\> ,$ contradicting the assumption
$U\supset K\cap S_0.$
\begin{pro}[Partial regularity in higher dimensions]\label{partial}
Let $u$ be a solution of (\ref{eq}) in $\Omega$ 
that is non-degenerate of second order. Then
the Hausdorff dimension of the set $S = \Omega\cap\{ u=0\}\cap \{\nabla u=0\}$ is less than or equal to $n-2\> .$
\end{pro}
\proof
Suppose that $s>n-2$ and that ${\cal H}^s(S) >0\> .$ Then we may
use \cite[Proposition 11.3]{giu}, Lemma \ref{dist}
as well as \cite[Lemma 11.5]{giu}
at ${\cal H}^s$-a.e. point of $S$ to obtain a blow-up limit
$u_0$ with the properties mentioned in Proposition \ref{fixedcenter},
satisfying ${\cal H}^{s,\infty}(S_0)>0$ 
for $S_0 := 
\{ u_0=0\}\cap \{\nabla u_0 = 0\}.$
According to Proposition \ref{fixedcenter} there are two possibilities:\\
1) $\Phi_{x^0}(0+)=-\infty$ and
$u_0$ is a homogeneous harmonic polynomial of degree $2.$
But for such a polynomial ${\cal H}^{n-2}(\{ u_0=0\}\cap \{\nabla u_0 = 0\})<+\infty$
and we obtain a contradiction. Thus the second possibility has to apply:\\
2) for some $\alpha\in (0,+\infty),$
$\alpha u_0$ is a solution of (\ref{eq}) on $\R^n$
that is homogeneous of degree $2.$ In this case we proceed
with the dimension reduction:
By \cite[Lemma 11.2]{giu}
we find a point $\bar x \in S_0 -\{0\}$ at which the density
in \cite[Proposition 11.3]{giu} is estimated from below.
Now each blow-up limit $u_{00}$ with respect to $\bar x$
(and with respect to a subsequence $m\to \infty$ such that 
the limit superior in \cite[Proposition 11.3]{giu} becomes a
limit) again satisfies the properties of Proposition \ref{fixedcenter};
in addition, we obtain from the homogeneity of $u_0$ as in
Lemma 3.1 of \cite{jga} that the rotated $u_{00}$ is constant
in the direction of the $n$-th unit vector. Defining 
$\bar u$ as the restriction of this rotated $u_{00}$ to
$\R^{n-1}\> ,$ it follows therefore that 
${\cal H}^{s-1}(\{ \bar u=0\}
\cap \{\nabla \bar u=0\})>0\> .$\\
Repeating the whole procedure $n-2$ times we obtain a nontrivial
homogeneous solution $u^\star$ of degree $2$ in $\R,$ satisfying
${\cal H}^{s-(n-2)}(\{ u^\star=0\}
\cap \{\nabla u^\star=0\})>0\> ,$ a contradiction.
\section{Lipschitz arcs in two dimensions}
In this section we show that the zero-set of the minimal
solution consists in every second-order non-degenerate
part of $\Omega$ of finitely many
Lipschitz arcs which end -- if so at all -- in quadruple junctions,
meeting at right angles.\\
In order to do the analysis, we have to prove uniform Lipschitz regularity close
to singular points. The difficulty is that convergence
to the blow-up limit is {\em not uniform} at singular
points. In \cite{SUW} we used a {\em novel} intersection-comparison
approach to obtain that close to the singular point the free boundary is {\em uniformly
the union of two graphs.}
In our case it turns out that the classical intersection-comparison
method (also called zero-number technique or lap-number technique) is sufficient,
when combined with a very elementary implicit function theorem argument.
The proof of the following theorem is inspired by
\cite{veronveron}.
\begin{theo}[Unique blow-up limit]\label{unique}
Let $n=2,$ let $u$ be the minimal solution of (\ref{eq}) in $\Omega$ 
and suppose that $u$ is
non-degenerate of second order at 
$x^0\in \Omega\cap \{ u=0\}\cap\{ \nabla u=0\}.$
Then, as $0<r\to 0,$
and $S(x^0,r) := \left(r^{1-n}\int_{\partial B_{r}(x^0)} u^2\> d{\cal H}^{n-1}
\right)^{1\over 2},$
$$u_r(x) := \frac{u(x^0+r x)}{S(x^0,r)}$$
converges to $p$ where $(p\circ U)(x)= (x_1^2-x_2^2)/\Vert x_1^2-x_2^2\Vert_{L^2(\partial B_1(0))}$ for some rotation
$U.$
\end{theo}
\proof
First, by Proposition \ref{fixedcenter},
for any $\tilde \epsilon>0$ there is $\tilde\rho>0$ such that
$$\dist(u_r,M_g) < \tilde \epsilon \hbox{ for }
r<\tilde \rho\; ;$$
here $M^*_g := \{ (x_1^2-x_2^2)/\Vert x_1^2-x_2^2\Vert_{L^2(\partial B_1(0))} \},
M_g := \{ q \> : \> q\circ V(x)= (x_1^2-x_2^2)/\Vert x_1^2-x_2^2\Vert_{L^2(\partial B_1(0))}$ 
for some rotation $V \},$
and $\dist(u_r,M_g) := \inf_{q\in M_g} \sup_{x\in B_1(0)} \vert u_r(x)-q(x)\vert.$
Denote by $U_\theta$ the counterclockwise rotation of positive angle
$\theta.$
If the statement of the theorem does not hold, then
-- by uniform
continuity of $t \mapsto u_{\exp(-t)}$ --
there exists a sequence $r_m\downarrow 0$ and
rotations $U_{\theta_1}$ and $U_{\theta_2}$ satisfying
$\vert\theta_{1}-\theta_{2}\vert = \epsilon\in (0,\pi)$
as well as
$$\dist(u_{r_{2m}}\circ U_{\theta_{1}}, M^*_g) \leq \tilde \epsilon
\hbox{ and } \dist(u_{r_{2m+1}}\circ U_{\theta_2}, M^*_g) \leq \tilde \epsilon$$
for $m=0,1,2,\dots$\\
Note that we may assume $\theta_1>\theta_2.$
Now let
$U = U_{{\theta_1+\theta_2\over 2}}, \omega= (\theta_1-\theta_2)/2\in (0,\pi/2)$ and define
$$
 \phi (r,\theta):=\frac{u(x^0 + rU(\cos \theta,\sin
 \theta))}{S(x^0,r)}\; .
 $$
For each $0<r<r_0$, the function $\phi(r,\cdot
 )$ defines a function on the unit circle $[-\pi,\pi).$
Inspired by applications of the Aleksandrov reflection (see for example
\cite{chen}, \cite{matanoveron}, \cite{veronveron})
we consider now
$$
 \xi (r,-\theta):=\phi(r,\theta) - \phi(r,-\theta).
 $$
Observe that $\xi(r,0)=\xi(r,\pi)=0.$
In what follows we will
prove that 
${\partial\xi\over {\partial\theta}}(r_{2m},0) \; < \; 0$
and
${\partial\xi\over {\partial\theta}}(r_{2m+1},0) \; > \; 0$
for large $m.$
The comparison principle (applied to the minimal solutions
$S(x^0,r)\phi(r,\theta)$ and $S(x^0,r)\phi(r,-\theta)$ in the
two-dimensional domain $[0,r_0)\times (0,\pi)$ with respect to the
original coordinates $x_1$ and $x_2$),
tells us now that the connected component of $\{ \xi <0\}$ touching
$(r_{2m},0)$ intersects $\{ r_0\}\times (0,\pi),$ and
that the connected component of $\{ \xi >0\}$ touching
$(r_{2m+1},0)$ intersects $\{ r_0\}\times (0,\pi).$
It follows that $\{ r_0\}\times (0,\pi)$ contains infinitely
many connected components of $\{ \xi >0\}$ and $\{ \xi <0\}.$
On the other hand we know that, provided that $r_0$ has been chosen
small enough,
$u_{r_0}\circ U_{\theta_1}$ is close to
$(x_1^2-x_2^2)/\Vert x_1^2-x_2^2\Vert_{L^2(\partial B_1(0))},$
so $\xi(r_0,\theta)$ is close to
$\zeta(\theta):=c_1(\cos(\omega+\theta)^2-\sin(\omega+\theta)^2 
-\cos(\omega-\theta)^2+\sin(\omega-\theta)^2)
= 2c_1(\cos(\theta+\omega)^2-\cos(\theta-\omega)^2)$
in $C^1([0,\pi]).$ As the zeroes of
$\zeta$ are all non-degenerate, it is not
possible that $\xi(r_0,\cdot)$ has infinitely many zeroes.
\\
Therefore, in order to finish the proof, we have to show
that ${\partial\xi\over {\partial\theta}}(r_{2m},0) \; < \; 0$
and
${\partial\xi\over {\partial\theta}}(r_{2m+1},0) \; > \; 0$
for large $m.$
As $\xi(r_{2m},\cdot)$ is close to
$\zeta$
in $C^1([0,\pi])$ and $\xi(r_{2m+1},\cdot)$ is close to
$-\zeta$
in $C^1([0,\pi])$
we need only calculate
$\zeta'(0)=-8c_1\cos(\omega)\sin(\omega)\le c_2(\omega)<0.$\\
So for $\tilde \epsilon$ and $r_0$ sufficiently small (depending on $\epsilon$)
we obtain a contradiction.
\begin{cor}[Lipschitz arcs]\label{arc}
Let $n=2,$ let $u$ be the minimal solution of (\ref{eq}) in $\Omega$ 
and suppose that $u$ is
non-degenerate of second order at
$x^0\in \Omega\cap \{ u=0\}\cap\{ \nabla u=0\}.$
Then $\{ u=0\}$ consists in an open neighborhood of $x^0$
of four Lipschitz arcs
meeting at right angles.
\end{cor}
\proof
By Theorem \ref{unique}
we know that for small $\delta>0$ and $r\in (0,\delta),$
the $u_r$ of Theorem \ref{unique}
satisfies
$$|\partial_2 (u_r\circ U)| \ge c_1 \hbox{ in } B_2(0)\cap
\{ |x_2| \ge 1/8\}$$
$$\hbox{and } B_2(0)\cap \{ u_r\circ U=0\}
\subset
B_{1/16}(\{ |x_1|=|x_2| \})\; .$$
Thus $(B_2(0)-B_{1/8}(0))\cap \{ u_r\circ U=0\}$ is the union of four
$C^1$-graphs $g_r^j (j=1,2,3,4)$ in the $x_2$-direction,
satisfying
$$ \Vert (g_r^j)'\Vert_{C^0(I_j)}
\le {1\over {c_1}}\Vert \nabla u_r\Vert_{C^0(B_2(0))}\le C_2\; ,$$
where $I_j=[\pm 1/2,\pm 1)$.
Rescaling yields the statement of the corollary.

%%%%%%%%%%%%%%%%%%%%%%%%%%%%%%%%%%%%%%%%%%%%%%%%%%%%%%%%%
\section{Regularity of Local Minimizers}

\begin{theo}[Regularity of local minimizers]\label{main}
Let $u$ be a minimizer of the energy
\[ E(v) = \int_{B_{r_0}(x^0)} |\nabla v|^2 - 2\max(v,0)\]
in $K:= \{ v\in W^{1,2}(B_{r_0}(x^0)) : v=u \hbox{ on } \partial
B_{r_0}(x^0)\}$.
Then the free boundary $\partial\left\{u>0\right\}$ is 
locally in $\Omega$ an analytic surface
and
$u|_{\overline{\{ u>0\}}}, u|_{\overline{\{ u<0\}}}$ are locally in $\Omega$
analytic functions.
\end{theo}
\begin{rem}
By Lemma \ref{maxmin} this implies the same regularity 
for the maximal solution.
\end{rem}
\begin{rem}
The theorem also implies that local minimizers are locally
in $\Omega$ of class $C^{1,1}$. Usually regularity of
the solution is proved
{\em before} proving regularity of the free boundary,
but here we do it the other way around.
\end{rem}
We start with some preliminary results.\\
\begin{lem}\label{lem:E}
Let $n\ge 2$, let $u$ be a minimizer of the energy
\[ E(v) = \int_{B_{r_0}(x^0)} |\nabla v|^2 - 2\max(v,0)\]
in $K:= \{ v\in W^{1,2}(B_{r_0}(x^0)) : v=u \hbox{ on } \partial
B_{r_0}(x^0)\}$, 
and suppose that
$\nabla u(x)\not=0$ on
$\left\{u=0\right\}\cap (B_{r_0}(x^0) \backslash \left\{x^0\right\})$. 
Then for $w\in C^\infty_0(B_{r_0}(x^0)\backslash \overline{B_\delta(x^0)})$ we have
$$0\le \delta^2 E(u)\cdot (w)(w) =  \int_{B_{r_0}(x^0)} |\nabla
w|^2 - \int_{\left\{u=0\right\}\cap
  B_{r_0}(x^0)}\frac{1}{|\nabla u|}w^2\; .$$
\end{lem}
\proof 
We define 
$$E_\epsilon(u)=\int_{B_{r_0}(x^0)}\frac12 |\nabla u|^2 - \gamma_\epsilon(u)$$
where $\gamma_\epsilon(u)$ is an approximation of $\max(u,0)$ such that
$\gamma_\epsilon''(u)=1/\epsilon$ if $u\in (0,\epsilon)$ and
zero otherwise.
We see that for $w\in C^\infty_0(B_{r_0}(x^0)\backslash \overline{B_\delta(x^0)})$ we have
$$\frac{1}{t^2}\left(E_\epsilon(u+tw) - E_\epsilon(u)-
  t\delta E_\epsilon(u)(w)\right)  = A^t_\epsilon$$
where
$$A^t_\epsilon=\frac{1}{t^2}\int_{B_{r_0}(x^0)} \frac{t^2}{2} |\nabla w|^2
-\left(\gamma_\epsilon(u+tw)-\gamma_\epsilon(u) -
  t\gamma_\epsilon'(u)w\right)$$
can be rewritten as
$$A^t_\epsilon= \int_{B_{r_0}(x^0)} \left(\frac{1}{2} |\nabla w|^2 -
\int_0^1d\alpha\ \int_0^\alpha d\sigma\ \gamma_\epsilon''(u+\sigma tw)\
w^2\right)$$
$$= \int_{B_{r_0}(x^0)} \frac{1}{2} |\nabla w|^2 -
\int_0^1d\alpha\ \int_0^\alpha {1\over \epsilon}
\int_{B_{r_0}(x^0)\cap \{0<u+\sigma tw<\epsilon\}}w^2\; .$$
By the co-area formula, we obtain for small $t$ that
$$A^t_\epsilon\to A^t_0=\int_{B_{r_0}(x^0)} \frac{1}{2} |\nabla w|^2 - 
  \int_0^1d\alpha\ \int_0^\alpha d\sigma\
  \int_{\left\{u+\sigma tw=0\right\}\cap B_{r_0}(x^0)}\frac{1}{|\nabla
  (u+\sigma tw)|}w^2$$
as $\epsilon\to 0$.
We conclude that 
$$0\le \frac{1}{t^2}\left(E(u+tw) - E(u)-
  t\delta E(u) (w)\right)  = A^t_0\; .$$
Last, we take the limit $t\to 0$ and obtain
$$0\le \frac12 \delta^2 E(u)(w)(w) = \int_{B_{r_0}(x^0)} \frac{1}{2} |\nabla
w|^2 - \frac{1}{2}\int_{\left\{u=0\right\}\cap
  B_{r_0}(x^0)}\frac{1}{|\nabla u|}w^2\; .$$
\begin{lem}
If $u$ is a solution in $B_{r_0}(0)$ of 
satisfying $\nabla u(0)=0$, then there exists a constant $C<\infty$ such that 
\begin{equation}\label{eq:chemin}
|\nabla u(x)|\le C |x|\log 1/|x| \; \hbox{ in } B_{r_0/2}(0)\; .  
\end{equation}
\end{lem}
\proof
From Chemin \cite{Chemin} we infer that
$\Delta u\in L^\infty$ implies $u\in C^2_*$. Thus
$\nabla u \in C^1_*$, and there exists a constant $C<\infty$ such that
$$|\nabla u(x)-\nabla u(y)|\le C|x-y| (1+\log 1/|x-y|)\; \hbox{ in } B_{r_0/2}(0)$$
which proves the Lemma.\\
{\sl Proof of Theorem \ref{main}:}\\
{\sl Step 1 (Dimension $n=2$):}\\
From Lemma \ref{minndeg} and Corollary \ref{isolated}
we know that locally the free boundary is either a $C^1$-arc
-- in which case the gradient is non-zero on the free boundary --
or a cross composed of $4$ Lipschitz arcs meeting at right angles 
-- in which case the gradient is zero at the center $x^1$ of the cross. 
We want to show that such a  cross is impossible for local minimizers.
We may assume $x^1=0$.\\
From Lemma \ref{lem:E} and (\ref{eq:chemin}), we deduce that for some
constant $c_1>0$
\begin{equation}\label{eq:eqchemin}
0\le \int_{B_{r_1}(0)} |\nabla
w|^2 - c_1\int_{\left\{u=0\right\}\cap
  B_{r_1}(0)}\frac{1}{|x|\log 1/|x|}w^2 \hbox{ for } w \in C^\infty_0(B_{r_1}(0)-B_{\delta r_1}(0))\; .  
\end{equation}
We now consider $w_\delta(x)=\phi(x)-\phi(\frac{x}{\delta})$
where $\phi\in C^\infty_0(B_{r_1}(0))$ such that $\phi=1$ on $B_{\frac{r_1}2}(0)$. It follows
that $\int_{B_{r_1}(0)} |\nabla w_\delta|^2 \le C_2<\infty$,
and using for large $i\in \N$ the regularity of $\{ u=0\}\cap (B_{2^{-i}}(0)-B_{2^{-i-1}}(0))$
as well as the closeness to the rotated cross we obtain
$$C_3 \int_{\left\{u=0\right\}\cap
  B_{r_1}(0)}\frac{1}{|x|\log 1/|x|}w_\delta^2 \ge 
-\int_{\delta}^{\frac{r_1}2} dr\ \frac{1}{r\log r} \to +\infty $$
 as $\delta\to 0$, a contradiction to the boundedness of $w_\delta$.
\\
Therefore the cross is not a local minimizer, and in dimension
$n=2$ the free boundary is locally in $\Omega$ a $C^1$-arc for each local minimizer.\\
\noindent {\sl Step 2 (Dimension $n>2$):}\\
We proceed by induction.\\
We assume that we have proved that for local minimizers, 
the free boundary is smooth up to the dimension $n-1\ge 2$.\\
Now we cannot have an accumulation of singularities in dimension $n$:
Blowing up at a limit point and blowing up a second time
at a singularity $\bar x\ne 0$ of the blow-up limit,
we would obtain as in the proof of 
Proposition \ref{partial} (see also \cite[Lemma 3.1, Lemma 3.2]{jga}) 
a local minimizer with a singularity in dimension $n-1$.
Thus singularities are isolated, and every blow-up limit of $u$
at each singularity is a harmonic polynomial of degree 2
whose gradient vanishes only at one point of the $0$-level set.
Thus we still have (\ref{eq:eqchemin}), and the free boundary is 
for large $i\in \N$ on $B_{2^{-i}}(0)-B_{2^{-i-1}}(0)$
close
to the zero level set of a homogeneous harmonic polynomial $P_i$ of degree
2 satisfying by Proposition \ref{fixedcenter}  
$|P_i(x)|\le C_4|x|^2$ and
$|\nabla P_i|\ge c_5 |x|$ on $(\R^n-\{ 0\})\cap \{ P_i=0\}$
where $C_4<\infty$ and $c_5>0$ do not depend on $i$.\\
Now we choose $w_\delta(x)=|x|^{-\left(\frac{n-2}{2}\right)}(\phi(x)-\phi(\frac{x}{\delta}))$
where $\phi\in C^\infty_0(B_{r_1}(0))$ and $\phi=1$ on $B_{\frac{r_1}2}(0)$.
We obtain
 $\int_{B_{r_1}(0)} |\nabla w_\delta|^2 \le C_6$, 
and using the regularity of $\{ u=0\}\cap (B_{2^{-i}}(0)-B_{2^{-i-1}}(0))$
as well as the closeness to $\{ P_i=0\}$
we see that
$$C_7\int_{\left\{u=0\right\}\cap
  B_{r_1}(0)}\frac{1}{|x|\log 1/|x|}w_\delta^2 \ge 
-\int_{\delta}^{\frac{r_1}2} dr\ \frac{1}{r\log r} \left(\inf_i \int_{\partial B_1(0)\cap
  \left\{P_i=0\right\}} 1 \> d{\cal H}^{n-2}\right) \to +\infty $$
as $\delta \to 0$, a contradiction to the boundedness of $w_\delta$.\\
Therefore the local minimizers have no singularities in dimension $n$, and
the free boundary is locally in $\Omega$ a $C^1$-surface.\\
\noindent {\em Step 3 : Analyticity of the free boundary}\\
We obtain analyticity of the free boundary as well as
analyticity of
$u|_{\overline{\{ u>0\}}},u|_{\overline{\{ u<0\}}}$
as in the proof of Theorem 4.1 in Chapter 6 of \cite{KS}.
See also \cite[Theorem 3.1']{KNS}.

%%%%%%%%%%%%%%%%%%%%%%%%%%%%%%%%%%%%%%%%%%%%%%%%%%%%%%%%%%%%%%%%%%%%%%%%%%%%%%%%

%%%%%%%%%%%%%%%%%%%%%%%%%%%%%%%%%%%%%%%%%%%%%%%%%%%%%%%%%
\section{The Cross Singularity}\label{counterex}
The reader may wonder whether there exists an example of a singularity
for the maximal solution (and thereby a counter-example
to the $W^{2,\infty}$-regularity in this unstable problem).
We have at this moment no conclusive answer, but
the following formal asymptotic expansion suggests
that the cross may be a possibility:

\begin{lem}\label{lem:imp}{\bf (Formal asymptotics)}\\
Let us assume the existence of a solution $u$ in $B_1(0)$
such that $u$ satisfies for
$$x_1=r\cos
\alpha,x_2=r\sin\alpha$$
$$u(-x_1,x_2)=u(x_1,x_2)=u(x_1,-x_2),$$
and the free boundary in the set $B_1(0)\cap \left\{x_1>0, x_2>0\right\}$ is given by
$$\alpha = \frac{\pi}{4} (1+\phi(\rho)),\quad \rho =
\frac{1}{-\log r}\; .$$
Moreover we assume that 
$$u>0\quad \mbox{in}\quad 0<\alpha <\frac{\pi}{4} (1+\phi(\rho))\quad
\mbox{and}\quad u<0 \quad \mbox{in}\quad \frac{\pi}{4}  (1+\phi(\rho)) <
\alpha <\pi/2\ .$$
Then formally $\phi$ and $u$ satisfy
$$\phi(\rho)=-\frac{\rho}{2} + O(\rho^2)$$
$$\hbox{and }u(x)=\frac{1}{2\pi}(x_1^2-x_2^2)(-\log |x|) +O(|x|^2)\ .$$
\end{lem}
\noindent {\sl Formal proof:}\\
We set
$$u(x_1,x_2)=\left\{\begin{array}{rll}
r^2(-\log(r))\ z^+(\rho,\theta^+) & \quad \mbox{for}\quad \displaystyle{\theta^+ =
\frac{\alpha}{1+\phi(\rho)}} & \quad
\mbox{if}\quad  \displaystyle{0<\alpha < \frac{\pi}{4}  (1+\phi(\rho))}\\
\\
-r^2(-\log(r))\ z^-(\rho,\theta^-) & \quad \mbox{for}\quad \displaystyle{\theta^- =
\frac{\pi/2-\alpha}{1-\phi(\rho)}} & \quad
\mbox{if}\quad  \displaystyle{\frac{\pi}{4}  (1+\phi(\rho))< \alpha <
  \pi/2} \ .
 \end{array} \right.$$
By continuity of $u$ and $\nabla u$ on the free boundary, we have
$$z^+=z^-=0\quad\mbox{and}\quad 
z^+_\theta=\left(\frac{1+\phi(\rho)}{1-\phi(\rho)}\right)z^-_\theta
\quad \mbox{for}\quad \theta^+=\theta^- = \pi/4\ .$$
By the symmetries of $u$ we also have 
$$z^+_\theta = z^-_\theta =0 \quad \mbox{for}\quad \theta^+=\theta^- =0\ .$$
Moreover
$z=z^+$ satisfies for $0< \theta=\theta^+ <\pi/4$:
$$z_{\theta\theta} + 4(1+\phi(\rho))^2 z + \left(\rho -4\rho
  z\right)(1+\phi(\rho))^2 +\rho^2 I[z,\phi] =0$$
where
$$\begin{array}{ll}
I[z,\phi]= & \displaystyle{4\left((1+\phi(\rho))^2 z_\rho -\theta
      z_\theta \phi'(\rho)(1+\phi(\rho))\right)}\\
& +\quad  \displaystyle{\rho^2\left((1+\phi(\rho))^2z_{\rho\rho} +
  \phi'^2(\rho)(\theta^2z_{\theta\theta} + 2 \theta z_\theta)
  -(1+\phi(\rho))\left\{2\theta z_{\rho\theta}\phi'(\rho)+\theta
    z_\theta\phi''(\rho)\right\}\right)}\; ,
\end{array}$$
and for $0< \theta=\theta^- <\pi/4$, the function
$z=z^-$ satisfies a similar equation with $\phi$ replaced by $-\phi$ and
without the term $\rho(1+\phi(\rho))^2$, i.e.
$$z_{\theta\theta} + 4(1-\phi(\rho))^2 z -4\rho
  z (1-\phi(\rho))^2 +\rho^2 I[z,-\phi] =0\ .$$
Let us now introduce the formal asymptotic expansion
$$\left\{\begin{array}{l}
\phi(\rho)=\phi^0+\rho\phi^1+\rho^2\phi^2+...\\
   \\
z^\pm=z^{\pm,0}(\theta)+\rho z^{\pm,1}(\theta)+\rho^2 z^{\pm,2}(\theta)+...
\end{array}\right.$$
where
$$\phi^0=0,\quad z^{\pm,0}(\theta)=A^0\cos(2\theta)\quad
\mbox{and}\quad A^0>0  \quad
\mbox{by assumption.}$$
For the order $0$ terms we obtain (for $0<\theta <\pi/4$)
$$z^{\pm,0}_{\theta\theta} +4z^{\pm,0}=0,\quad z^{\pm,0}_\theta(0)=0, \quad z^{\pm,0}(\pi/4)=0, \quad
z^{+,0}_\theta(\pi/4)=z^{-,0}_\theta(\pi/4)\ ,$$
which is compatible with the assumptions.
For the order $1$ terms we obtain (for $0<\theta <\pi/4$)
$$\left\{\begin{array}{l}
z^{+,1}_{\theta\theta} +4z^{+,1} + (8\phi^1-4)z^{+,0} +1 =0,\quad \quad 
z^{-,1}_{\theta\theta} +4z^{-,1} + (-8\phi^1-4)z^{+,0} =0,\\
\\
z^{\pm,1}_\theta(0)=0, \quad z^{\pm,1}(\pi/4)=0, \quad 
z^{+,1}_{\theta}(\pi/4)=z^{-,1}_{\theta}(\pi/4) -2\phi^1
z^{-,0}_{\theta}(\pi/4)\ .   \end{array}\right.$$
Thus
$$\begin{array}{l}
\displaystyle{z^{+,1}(\theta)=A^{+,1}\cos(2\theta) +(1-2\phi^1)A^0\theta\sin(2\theta)-
\frac14}\\
\displaystyle{\hbox{and }z^{-,1}(\theta)=A^{-,1}\cos(2\theta) + (1+2\phi^1)A^0\theta\sin(2\theta)}
\end{array}$$
with $\phi^1=-\frac12$, $A^0=\frac{1}{2\pi}$ and
$A^{+,1}-A^{-,1}=\frac{1}{\pi}$.\\

%%%%%%%%%%%%%%%%%%%%%%%%%%%%%%%%%%%%%%%%%%%%%%%%%%%%%%%%%
\section{Open Questions}
The most urgent remaining questions in this context
are, whether the cross singularity can be proven to exist,
and whether there are examples of second order degeneracy.
\\
Another interesting point is whether methods similar
to those used in this paper (possibly combined with
arguments as in \cite{aguilera}) can be used to
prove regularity in the composite membrane problem
(see \cite{blank}, \cite{chanillo1}, \cite{chanillo2}).  
%%%%%%%%%%%%%%%%%%%%%%%%%%%%%%%%%%%%%%%%%%%%%%%%%%%%%%%%%%%%%%%%%%%%%%%%%%%%%%%%

\bibliographystyle{plain}
\bibliography{monneauweiss050715.bib}

\def\cprime{$'$} \def\cprime{$'$}
\begin{thebibliography}{10}

\bibitem{aguilera}
N.~Aguilera, H.~W. Alt, and L.~A. Caffarelli.
\newblock An optimization problem with volume constraint.
\newblock {\em SIAM J. Control Optim.}, 24(2):191--198, 1986.

\bibitem{khaikin}
A.P. Aldushin and B.I. Khaikin.
\newblock Combustion of mixtures forming condensed reaction-products.
\newblock {\em Combus. Explos. Shock Waves}, 10:273--280, 1974.

\bibitem{ACF}
Hans~Wilhelm Alt, Luis~A. Caffarelli, and Avner Friedman.
\newblock Variational problems with two phases and their free boundaries.
\newblock {\em Trans. Amer. Math. Soc.}, 282(2):431--461, 1984.

\bibitem{beck}
J.~M. Beck and V.~A. Volpert.
\newblock Nonlinear dynamics in a simple model of solid flame microstructure.
\newblock {\em Phys. D}, 182(1-2):86--102, 2003.

\bibitem{veronveron}
Marie-Fran{\c{c}}oise Bidaut-V{\'e}ron, Victor Galaktionov, Philippe Grillot,
  and Laurent V{\'e}ron.
\newblock Singularities for a 2-dimensional semilinear elliptic equation with a
  non-{L}ipschitz nonlinearity.
\newblock {\em J. Differential Equations}, 154(2):318--338, 1999.

\bibitem{blank}
Ivan Blank.
\newblock Eliminating mixed asymptotics in obstacle type free boundary
  problems.
\newblock {\em Comm. Partial Differential Equations}, 29(7-8):1167--1186, 2004.

\bibitem{stars}
Luis~A. Caffarelli and Avner Friedman.
\newblock The shape of axisymmetric rotating fluid.
\newblock {\em J. Funct. Anal.}, 35(1):109--142, 1980.

\bibitem{cazenavelions}
Thierry Cazenave and Pierre-Louis Lions.
\newblock Solutions globales d'\'equations de la chaleur semi lin\'eaires.
\newblock {\em Comm. Partial Differential Equations}, 9(10):955--978, 1984.

\bibitem{chanillo2}
S.~Chanillo, D.~Grieser, M.~Imai, K.~Kurata, and I.~Ohnishi.
\newblock Symmetry breaking and other phenomena in the optimization of
  eigenvalues for composite membranes.
\newblock {\em Comm. Math. Phys.}, 214(2):315--337, 2000.

\bibitem{chanillo1}
S.~Chanillo, D.~Grieser, and K.~Kurata.
\newblock The free boundary problem in the optimization of composite membranes.
\newblock In {\em Differential geometric methods in the control of partial
  differential equations (Boulder, CO, 1999)}, volume 268 of {\em Contemp.
  Math.}, pages 61--81. Amer. Math. Soc., Providence, RI, 2000.

\bibitem{Chemin}
Jean-Yves Chemin.
\newblock {\em Perfect incompressible fluids}, volume~14 of {\em Oxford Lecture
  Series in Mathematics and its Applications}.
\newblock The Clarendon Press Oxford University Press, New York, 1998.

\bibitem{chen}
Xu-Yan Chen.
\newblock Uniqueness of the {$\omega$}-limit point of solutions of a semilinear
  heat equation on the circle.
\newblock {\em Proc. Japan Acad. Ser. A Math. Sci.}, 62(9):335--337, 1986.

\bibitem{matanoveron}
Xu-Yan Chen, Hiroshi Matano, and Laurent V{\'e}ron.
\newblock Anisotropic singularities of solutions of nonlinear elliptic
  equations in {${\bf R}\sp 2$}.
\newblock {\em J. Funct. Anal.}, 83(1):50--97, 1989.

\bibitem{giu}
Enrico Giusti.
\newblock {\em Minimal surfaces and functions of bounded variation}, volume~80
  of {\em Monographs in Mathematics}.
\newblock Birkh\"auser Verlag, Basel, 1984.

\bibitem{KNS}
D.~Kinderlehrer, L.~Nirenberg, and J.~Spruck.
\newblock Regularity in elliptic free boundary problems.
\newblock {\em J. Analyse Math.}, 34:86--119 (1979), 1978.

\bibitem{KS}
David Kinderlehrer and Guido Stampacchia.
\newblock {\em An introduction to variational inequalities and their
  applications}, volume~88 of {\em Pure and Applied Mathematics}.
\newblock Academic Press Inc. [Harcourt Brace Jovanovich Publishers], New York,
  1980.

\bibitem{solomonov}
V.A. Knyazik, A.G. Merzhanov, V.B. Solomonov, and A.S. Shteinberg.
\newblock Macrokinetics of high-temperature titanium interaction with carbon
  under electrothermal explosion conditions.
\newblock {\em Combus. Explos. Shock Waves}, 21:333--337, 1985.

\bibitem{henrik}
Henrik Shahgholian.
\newblock The singular set for the composite membrane problem.
\newblock {\em Preprint}.

\bibitem{SUW}
Henrik Shahgholian, Nina Uraltseva, and Georg~S. Weiss.
\newblock Global solutions of an obstacle-problem-like equation with two
  phases.
\newblock {\em Monatsh. Math.}, 142(1-2):27--34, 2004.

\bibitem{varma}
A.~Varma, A.S. Mukasyan, and S.~Hwang.
\newblock Dynamics of self-propagating reaction in heterogeneous media:
  experiments and model.
\newblock {\em Chem. Eng. Sci.}, 56:1459--1466, 2001.

\bibitem{interfaces}
G.~S. Weiss.
\newblock An obstacle-problem-like equation with two phases: pointwise
  regularity of the solution and an estimate of the {H}ausdorff dimension of
  the free boundary.
\newblock {\em Interfaces Free Bound.}, 3(2):121--128, 2001.

\bibitem{cpde}
Georg~S. Weiss.
\newblock Partial regularity for weak solutions of an elliptic free boundary
  problem.
\newblock {\em Comm. Partial Differential Equations}, 23(3-4):439--455, 1998.

\bibitem{jga}
Georg~Sebastian Weiss.
\newblock Partial regularity for a minimum problem with free boundary.
\newblock {\em J. Geom. Anal.}, 9(2):317--326, 1999.

\bibitem{zik}
O.~Zik, Z.~Olami, and E.~Moses.
\newblock {\em Phys. Rev. Lett.}, 81:3868 pp., 1998.

\end{thebibliography}

\end{document}